\documentclass[12pt]{amsart}
\usepackage{amsmath,amscd,amssymb,amsfonts,graphics}
\usepackage{hyperref}
\newcommand{\hl}{\hyperlink}
\newcommand{\htt}{\hypertarget}
\newcommand{\h}{\hbox}
\newcommand{\q}{\quad}

\newcommand{\bs}{\par\bigskip}
\newcommand{\ms}{\par\medskip}
\newcommand{\sk}{\par\smallskip}
\newcommand{\bsn}{\par\bigskip\noindent}
\newcommand{\msn}{\par\medskip\noindent}
\newcommand{\skn}{\par\smallskip\noindent}
\newcommand{\ges}{\geqslant}
\newcommand{\les}{\leqslant}
\newcommand{\1}{\hskip1pt}

\newcommand{\msum}{\hbox{$\sum$}}
\newcommand{\mopl}{\hbox{$\bigoplus$}}

\newcommand{\B}{{\mathcal B}}
\newcommand{\D}{{\mathcal D}}
\newcommand{\F}{{\mathcal F}}
\newcommand{\G}{{\mathcal G}}
\newcommand{\J}{{\mathcal J}}
\newcommand{\I}{{\mathcal I}}
\newcommand{\K}{{\mathcal K}}
\newcommand{\Hc}{{\mathcal H}}

\newcommand{\M}{{\mathcal M}}
\newcommand{\OO}{{\mathcal O}}
\newcommand{\Sc}{{\mathcal S}}

\newcommand{\Q}{{\mathbb Q}}
\newcommand{\C}{{\mathbb C}}

\newcommand{\R}{{\mathbb R}}
\newcommand{\RR}{{\mathbf R}}
\newcommand{\DD}{{\mathbf D}}
\newcommand{\Z}{{\mathbb Z}}

\newcommand{\alt}{\widetilde{\alpha}}
\newcommand{\bt}{\widetilde{b}}
\newcommand{\Fti}{\widetilde{F}}
\newcommand{\Ht}{\widetilde{H}}
\newcommand{\Vt}{\widetilde{V}}
\newcommand{\Yt}{\widetilde{Y}}
\newcommand{\Xt}{\widetilde{X}}

\newcommand{\jt}{{}\,\widetilde{\!j}{}}

\newcommand{\al}{\alpha}

\newcommand{\om}{\omega}
\newcommand{\Om}{\Omega}
\newcommand{\Omu}{\underline{\Omega}}
\newcommand{\DR}{{\rm DR}}
\newcommand{\Diff}{{\rm Diff}}
\newcommand{\DbcF}{D^b_{\rm coh}F}
\newcommand{\codim}{{\rm codim}}
\newcommand{\Gr}{{\rm Gr}}
\newcommand{\Sp}{{\rm Sp}}
\newcommand{\Sing}{{\rm Sing}}
\newcommand{\dd}{\partial}
\newcommand{\ddd}{{\rm d}}
\newcommand{\mm}{{\mathfrak m}}
\newcommand{\bl}{\bigl}
\newcommand{\br}{\bigr}
\newcommand{\pl}{\1{+}\1}
\newcommand{\mi}{\1{-}\1}
\newcommand{\eq}{\,{=}\,}

\newcommand{\less}{\,{\leqslant}\,}
\newcommand{\gess}{\,{\geqslant}\,}

\newcommand{\ssb}{\raise.15ex\h{${\scriptscriptstyle\bullet}$}}
\newcommand{\ssc}{\,\raise.15ex\h{${\scriptstyle\circ}$}\,}
\newcommand{\onto}{\twoheadrightarrow}
\newcommand{\into}{\hookrightarrow}
\newcommand{\simto}{\,\,\rlap{\hskip1.5mm\raise1.4mm\hbox{$\sim$}}\hbox{$\longrightarrow$}\,\,}
\newcommand{\plim}{\rlap{\raise-5.5pt\h{$\,\leftarrow$}}{\rm lim}}
\begin{document}
\title{Higher Du Bois singularities of hypersurfaces}
\author[S.-J. Jung]{Seung-Jo Jung}
\address{S.-J. Jung : Department of Mathematics Education, and Institute of Pure and Applied Mathematics, Jeonbuk National University, Jeonju, 54896, Korea}
\email{seungjo@jbnu.ac.kr}
\author[I.-K. Kim]{In-Kyun Kim}
\address{I.-K. Kim : Department of Mathematics, Yonsei University, Seoul, Korea}
\email{soulcraw@gmail.com}
\author[M. Saito]{Morihiko Saito}
\address{M. Saito : RIMS Kyoto University, Kyoto 606-8502 Japan}
\email{msaito@kurims.kyoto-u.ac.jp}
\author[Y. Yoon]{Youngho Yoon}
\address{Y. Yoon : Department of Mathematics, Chungnam National University, 99 Daehak-ro, Daejeon 34134, Korea}
\email{mathyyoon@gmail.com}
\thanks{This work was partially supported by NRF grant funded by the Korea government(MSIT) (the first author: NRF-2021R1C1C1004097, the second author: NRF-2020R1A2C4002510, and the fourth author: NRF-2020R1C1C1A01006782). The third author was partially supported by JSPS Kakenhi 15K04816.}
\subjclass{32S35}
\begin{abstract} For a complex algebraic variety $X$, we introduce higher $p$-Du Bois singularity by imposing canonical isomorphisms between the sheaves of K\"ahler differential forms $\Omega_X^q$ and the shifted graded pieces of the Du Bois complex $\underline{\Omega}_X^q$ for $q\le p$. If $X$ is a reduced hypersurface, we show that higher $p$-Du~Bois singularity coincides with higher $p$-log canonical singularity, generalizing a well-known theorem for $p=0$. The assertion that $p$-log canonicity implies $p$-Du Bois has been proved by Mustata, Olano, Popa, and Witaszek quite recently as a corollary of two theorems asserting that the sheaves of reflexive differential forms $\Omega_X^{[q]}$ ($q\le p$) coincide with $\Omega_X^q$ and $\underline{\Omega}_X^q$ respectively, and these are shown by calculating the depth of the latter two sheaves. We construct explicit isomorphisms between $\Omega_X^q$ and $\underline{\Omega}_X^q$ applying the acyclicity of a Koszul complex in a certain range. We also improve some non-vanishing assertion shown by them using mixed Hodge modules and the Tjurina subspectrum in the isolated singularity case. This is useful for instance to estimate the lower bound of the maximal root of the reduced Bernstein-Sato polynomial in the case where a quotient singularity is a hypersurface and its singular locus has codimension at most 4.
\end{abstract}
\maketitle
\centerline{\bf Introduction}
\bsn
Let $X$ be a reduced closed subvariety of a smooth complex algebraic variety $Y$. We have a bounded complex of mixed Hodge module denoted by $\Q_{h,X}$ in this paper and such that its underlying $\Q$-complex is the constant sheaf $\Q_{X^{\rm an}}$ and $H^0(X,\Q_{h,X})$ has weight 0, see for instance \cite[(4.4.2)]{mhm} and also a remark after \cite[Cor.\,1]{FPS}, etc.
\sk
We denote by $(M_X^{\prime\ssb},F)$ the underlying complex of filtered left $\D_Y$-modules of $\Q_{h,X}$. We have the associated de Rham complex $\DR_Y(M_X^{\prime\ssb},F)$ (which is shifted by $d_Y$ to the left as usual, where $d_Y:=\dim Y$). Let $\Omu_X^p$ be the $p$\1th graded piece of the Du~Bois complex shifted by $p$, see \cite{DB}. (Its construction is essentially the same as in \cite[\S 8.1]{De}.) As a corollary of \cite[Thm.\,4.2]{mhc} (see also Theorem~\hl{T1.5}{1.5} below), there are canonical isomorphisms for $p\in\Z$
\htt{1}{}
$$\Omu_X^p[-p]=\Gr_F^p\DR_Y(M_X^{\prime\ssb})\q\h{in}\,\,\,D^b_{\rm coh}(\OO_Y).
\leqno(1)$$
\sk
We say that $X$ has only {\it higher $\1 p$-Du~Bois singularities\1} if there are canonical isomorphisms
\htt{2}{}
$$\Om_X^q\,\simto\,\Omu_X^q\q\h{in}\,\,\,D^b_{\rm coh}(\OO_X)\q(\forall\,q\in[0,p]),
\leqno(2)$$
see also Remarks~\hl{R2.2}{2.2}--\hl{R2.3}{3} below. Here the $\Om_X^q$ are sheaves of K\"ahler differential forms, see 1.6 below. For $p\eq0$, this coincides with the definition of Du~Bois singularity. In the hypersurface case the above isomorphisms are unique as long as their restrictions to the smooth part of $X$ are canonical isomorphisms, since we can show that the $\Om_X^q$ are {\it torsion-free\1} for $q\in[0,p]$, see Proposition~\hl{P2.2}{2.2} below. Note that singularity becomes {\it milder\1} as $p$ increases, and we get {\it smoothness\1} for $p>d_Y/2$ in the hypersurface case by Theorem~\hl{T1}{1} together with (\hl{3}{3}) below.
\sk
From now on, assume $X\subset Y$ is a reduced hypersurface defined by a function $f$. In this case it is well-known that $M_X^{\prime\ssb}[d_X]$ is a mixed Hodge module, which will be denoted by $M_X$. Let $-\alt_X$ be the maximal root of the reduced Bernstein-Sato polynomial $\bt_f(s):=b_f(s)/(s{+}1)$. This depends only on $X$, see for instance \cite[Rem.\,4.2\,(i)]{wh}. Note that $\alt_X$ is a positive rational number if $X$ is singular, see \cite{Ka}. Set $\alt_X:=+\infty$ if $X$ is nonsingular. By \cite[Thm.\,0.4]{mic}, we have
\htt{3}{}
$$\alt_X\les d_Y/2\q\h{if}\q\alt_X\,{<}\,{+}\infty.
\leqno(3)$$
Apply this to the restriction of $X$ to a sufficiently general transversal slice to the smooth part of $\Sing\,X$. It is well-known that $\alt_X$ does not change under this restriction locally, see for instance \cite[Lem.\,2.2]{DMST} (using (\hl{12}{12}) below and GAGA). We thus get the following (which seems to be proved in \cite{MOPW} using a {\it slightly different\1} method).
\par\htt{P1}{}\msn
{\bf Proposition~1.} {\it In the above notation and assumption we have}
\htt{4}{}
$$\alt_X\les \tfrac{1}{2}\,\codim_Y\Sing\, X.
\leqno(4)$$
\ms
We say that $X$ has only (higher) $p$-{\it log canonical singularities\1} ($p\in\Z_{>0}$) if $I_p(X)\eq\OO_Y$, see \cite{MP}. Here the $I_p(X)$ are Hodge ideals. (These can be defined for the associated analytic spaces, see \cite{JKSY}). They are decreasing for $p$, since $X$ is a $\Z$-divisor. We have
\htt{5}{}
$$F_p\OO_Y(*X)=I_p(X)f^{-p-1}\q\h{in}\,\,\,\OO_Y(*X)\q(p\in\Z_{\ges 0}),
\leqno(5)$$
by definition, where the left-hand side is the Hodge filtration of the left $\D_Y$-module $\OO_Y(*X)$ underling the mixed Hodge module $(j_U)_*\Q_{h,U}$ with $U:=Y\,{\setminus}\,X\buildrel{j_U}\over\into Y$, see for instance \cite[Prop.\,2.11]{mhm}.
\sk
In this paper we prove the following.
\par\htt{T1}{}\msn
{\bf Theorem~1.} {\it For a reduced hypersurface $X\subset Y$, the following three conditions for $p\in\Z_{>0}$ are equivalent to each other\,$:$
\par\htt{a}{}\skn
\rlap{\rm(a)}\hskip.8cm\hangindent=.8cm\hangafter=1
$\alt_X\gess p{+}1$,
\par\htt{b}{}\skn
\rlap{\rm(b)}\hskip.8cm\hangindent=.8cm\hangafter=1
$X$ has only higher $p$-Du~Bois singularities,
\par\htt{c}{}\skn
\rlap{\rm(c)}\hskip.8cm\hangindent=.8cm\hangafter=1
$X$ has only higher $p$-log canonical singularities.}
\ms
The equivalence (\hl{a}{a})\,$\Leftrightarrow$\,(\hl{c}{c}) was proved in \cite[Cor.\,1]{hi} using GAGA. (This follows from (\hl{11}{11}--\hl{12}{12}) below.) The implication (\hl{a}{a}) $\Rightarrow$ (\hl{b}{b}) has been shown in \cite{MOPW} quite recently as a corollary of two theorems asserting that the sheaves of reflexive differential forms $\Om_X^{[q]}$ (see for instance \cite{Keb}, \cite{Mic}) for $q\le p$ coincide with $\Om_X^q$ and $\Omu_X^q$ respectively. These are proved by calculating the depth of the latter two sheaves. We construct explicit isomorphisms between $\Om_X^q$ and $\Omu_X^q$ applying the acyclicity of a Koszul complex in a certain range. More precisely, we prove the equivalence (\hl{b}{b})\,$\Leftrightarrow$\,(\hl{c}{c}) using (\hl{1}{1}) and Propositions~\hl{P1}{1}--\hl{P2}{2} as well as the compatibility of the de Rham functor $\DR_Y$ with the functor $\DD$ (here the equivalence (\hl{a}{a})\,$\Leftrightarrow$\,(\hl{c}{c}) is also employed to apply (\hl{3}{3})). For $p\eq0$, the equivalence (\hl{a}{a})\,$\Leftrightarrow$\,(\hl{b}{b}) was shown in \cite{mos}, see \cite{KS} for (\hl{b}{b})\,$\Leftrightarrow$\,(\hl{c}{c}) with $p\eq 0$. (Here it seems more correct to call the {\it pair\1} $(Y,X)$ log canonical, since an embedded resolution is used in the definition, see Remark~\hl{R2.3}{2.3} below for higher $p$-log canonicity.)
\sk
For the proof of Theorem~\hl{T1}{1}, the following is also needed.
\par\htt{P2}{}\msn
{\bf Proposition~2.} {\it We have the vanishing}
\htt{6}{}
$$\Hc^p(\Om_Y^{\ssb},\ddd f\wedge)=0\q\h{if}\q p<\codim_Y\Sing\,X.
\leqno(6)$$
\ms
Here $(\Om_Y^{\ssb},\ddd f\wedge)$ is locally the Koszul complex associated to the partial derivatives $\dd f/\dd y_i$ with $y_1,\dots,y_{d_Y}$ local coordinates of $Y$. See for instance \cite[Thm.\,16.8]{Mat} (and \cite[Rem.\,1.9\,(iv)]{DiSa1}). This vanishing is quite well-known in the hypersurface singularity theory, and is sometimes used even without any reference, see for instance \cite{DiSt}, \cite[(3.1.2)]{nwh}.
\sk
Theorem~\hl{T1}{1} might be viewed as a rather immediate consequence of \cite[Thm.\,4.2]{mhc} (see also Theorem~\hl{T1.5}{1.5} below) together with Propositions~\hl{P1}{1}--\hl{P2}{2} and (\hl{11}{11}--\hl{12}{12}), see 2.2--3 below.
\sk
In this paper we prove also the following (which seems quite implicit in the proof of \cite[Thm.\,1.5]{MOPW}).
\par\htt{T2}{}\msn
{\bf Theorem~2.} {\it If $\alt_X\,{>}\,p{+}1$ with $X$ a reduced hypersurface defined by a function $f$, there is an isomorphism in} $D^b_{\rm coh}(\OO_Y)\,{:}$
\htt{7}{}
$$\Omu_X^{d_X-p}[p]=\bl(\sigma_{\ges d_Y-p}(\Om_Y^{\ssb}|_X,\ddd f\wedge)\br)[d_Y].
\leqno(7)$$
\sk
Here $\sigma_{\ges k}$ is the ``stupid truncation", and $|_X$ means the tensor product with $\OO_X$ over $\OO_Y$. We deduce this theorem and Theorem~\hl{T3}{3} below from the canonical isomorphisms in Theorem~\hl{T3.2}{3.2} below. The latter theorem is shown using a calculation related to the nearby and vanishing cycle Hodge modules, see 3.1 below. Here \cite[Cor.\,1]{hi}, that is, (\hl{a}{a})\,$\Leftrightarrow$\,(\hl{c}{c}) in Theorem~\hl{T1}{1}, is also needed, and the assumption $\alt_X\,{>}\,p{+}1$ in Theorem~\hl{T2}{2} cannot be relaxed as in Theorem~\hl{T3}{3} below.
\sk
Set
$$\aligned&\K_{\ddd f}^p:={\rm Ker}(\ddd f{\wedge}\,{:}\,\Om_Y^p\,{\to}\,\Om_Y^{p+1}),\,\,\,\I_{\ddd f}^p:=\ddd f{\wedge}\,\Om_Y^{p-1},\,\,\,\Hc_{\ddd f}^p:=\K_{\ddd f}^p/\I_{\ddd f}^p,\\&{}_f\M:=\M/f\M,\q {}^f\!\M:={\rm Ker}(f:\M\to\M),\endaligned$$
for an $\OO_Y$-module $\M$ in general.
\sk
Since the tensor product with $\OO_X$ over $\OO_Y$ for a locally free $\OO_Y$-module is given by the mapping cone of multiplication by $f$, we get the following (using the spectral sequence associated to the double complex whose first differential is the multiplication by $f$ so that it degenerates at $E_2$).
\par\htt{C1}{}\msn
{\bf Corollary~1.} {\it If $\alt_X\,{>}\,p{+}1$ with $X$ as in Theorem~$\hl{T2}{2}$ above, there are short exact sequences
\htt{8}{}
$$0\to{}_f\K_{\ddd f}^{d_Y-p}\to\Hc^0(\Omu_X^{d_X-p})\to{}^f\Hc_{\ddd f}^{d_Y-p+1}\to0,
\leqno(8)$$
\vskip-6mm
\htt{9}{}
$$0\to{}_f\Hc_{\ddd f}^{d_Y-p+j}\to\Hc^j(\Omu_X^{d_X-p})\to{}^f\Hc_{\ddd f}^{d_Y-p+j+1}\to0,
\leqno(9)$$
for $j\in[1,p{-}1]$ together with an isomorphism
\htt{10}{}
$$\Hc^p(\Omu_X^{d_X-p})={}_f\Hc_{\ddd f}^{d_Y},
\leqno(10)$$
where $\Hc^j(\Omu_X^{d_X-p})=0$ for $j\notin[0,p]\,\,($and $\,d_Y\eq d_X{+}1)$.}
\ms
This gives an improvement of \cite[Thm.\,1.5]{MOPW}. Note that ${}_f\Hc_{\ddd f}^{d_Y}\cong\om_Y/((\dd f),f)\om_Y$, and $\bl((\dd f)\,{:}\,f\br)/(\dd f)={\rm Ker}\bl(f\,{:}\,\OO_Y/(\dd f)\,{\to}\,\OO_Y/(\dd f)\br)$, where $(\dd f)$ denotes the Jacobian ideal generated by the partial derivatives $\dd f/\dd y_i$.
In the isolated singularity case we have $\Hc_{\ddd f}^j\eq 0$ for $j\ne d_Y$. (This is a special case of Proposition~\hl{P2}{2}.)
\sk
Let $\Vt$ be the {\it microlocal\1} $V$-filtration on $\OO_Y$ associated to $f$. This is a decreasing filtration indexed by $\Q$, see \cite{mic}. Set $\,\Gr_{\Vt}^{\al}\OO_Y:=\Vt^{\al}\OO_Y/\Vt^{>\al}\OO_Y\,$ ($\al\in\Q$), and
$$\alt_f^{\,\rm min.int}:=\min\bl\{p\in\Z\,\big|\,\Gr_{\Vt}^p\OO_Y\ne 0\br\}.$$
This can be larger than the absolute value of the maximal integral root of $b_f(s)/(s{+}1)$, for instance, if $f\eq x^6{+}y^5{+}x^3y^3{+}z^5{+}w^3$. This phenomenon is closely related to well-known {\it integral shifts\1} between the spectral numbers and the roots of $b_f(s)/(s{+}1)$ in the isolated hypersurface singularity case (especially for semi-weighted-homogenous polynomials), see Remark~\hl{R3.4d}{3.4d} below.  By \cite[Thm.\,1 and (1.3.8)]{hi}, we have
\htt{11}{}
$$I_p(X)=\Vt^{p+1}\OO_Y\,\,\,\,\h{mod}\,\,\,\,(f)\q(p\in\Z),
\leqno(11)$$
\vskip-7mm
\htt{12}{}
$$\alt_X=\min\bl\{\al\in\Q\,\big|\,\Gr_{\Vt}^{\al}\OO_Y\ne 0\br\}\les\alt_f^{\,\rm min.int}.
\leqno(12)$$
\q
Using these, we can prove the following.
\par\htt{T3}{}\msn
{\bf Theorem~3.} {\it Let $X\subset Y$ be a reduced hypersurface defined by $f$. Assume $\alt_f^{\,\rm min.int}\,{>}\,p{+}1$, and $\Vt^{p+1}\OO_Y\not\subset((\dd f),f)$ in $\1\OO_Y$. Then $\Hc^p\Omu_X^{d_X-p}\,{\ne}\,0$.}
\ms
The first assumption is satisfied for any $p$ (that is, $\alt_f^{\,\rm min.int}={+}\infty$) if $X$ is a $\Q$-homology manifold, or equivalently, if the unipotent monodromy part of the vanishing cycle complex $\varphi_{f,1}\Q_{Y^{\rm an}}$ vanishes, see Theorem~\hl{T3.1}{3.1} below (and \cite{Mi} for the isolated singularity case).
\sk
In the case $X$ has only one singular point, let $\Sp_f(t)=\msum_{i=1}^{\mu_f}\,t^{\1\al_{f,i}}$ be the {\it spectrum\1} of $f$, where $\mu_f$ is the Milnor number, see \cite{St1}, \cite{JKSY}. The $\al_{f,i}$ are called the {\it spectral numbers\1} of $f$, see also Remark~\hl{R3.5}{3.5} below. We have the equality
\htt{13}{}
$$\alt_f^{\,\rm min.int}=\min\{\al_{f,i}\,|\,\al_{f,i}\,{\in}\,\Z\}.
\leqno(13)$$
The right-hand side is called the {\it minimal integral spectral number.} We have the {\it Tjurina subspectrum\1} $\Sp^{\rm Tj}_f(t)\eq\msum_{i\in J_f}\,t^{\al_{f,i}}$ where $J_f$ is a certain subset of $\{1,\dots,\mu_f\}$ with $|J_f|\eq\tau_f$ the Tjurina number of $f$, see \cite{JKSY}. The $\al_{f,i}$ for $i\in J_f$ are called the {\it Tjurina subspectral numbers.} We denote by $\al^{\rm max.Tj}_f$ the {\it maximal Tjurina subspectral number.}
\sk
By \cite[Prop.\,1.4]{JKSY} (see also a remark after (\hl{3.4.1}{3.4.1}) below), Theorem~\hl{T3}{3} implies the following.
\par\htt{C2}{}\msn
{\bf Corollary~2.} {\it Assume a hypersurface $X$ has only one singular point with $\alt_f^{\,\rm min.int}\,{>}\,p{+}1$, $\al^{\rm max.Tj}_f\gess p{+}1$. Then $\Hc^p\Omu_X^{d_X-p}\,{\ne}\,0$.}
\ms
From Corollary~\hl{C2}{2} we can deduce for instance the following.
\par\htt{C3}{}\msn
{\bf Corollary~3.} {\it If a quotient singularity $X$ is a hypersurface and $\codim_X\Sing\,X\less 4$, then $\alt_X\in\bl(1,\tfrac{3}{2}\br]$.}
\ms
Corollary~\hl{C3}{3} improves an estimate $\alt_X\in(1,2\1]$ proved in \cite[Cor.\,1.6]{MOPW}, if $X$ satisfies an additional condition $\codim_X\Sing\,X\less 4$ (which may be conjectured to be unnecessary). We can get the strict inequality $\alt_X\,{<}\,2$ from (\hl{10}{10}) together with \cite[Thm.\,5.3]{DB}, since $V$-manifolds are $\Q$-homology manifolds. The lower bound $\alt_X\,{>}\,1$ comes from \cite[Thm.\,0.4]{rat}, since it is well-known that quotient singularities are rational singularities. The upper bound $\alt_X\eq\tfrac{3}{2}$ can be attained by $f\eq x^2\pl y^2\pl z^2$. It does not seem easy to construct a quotient singularity $X$ which is a hypersurface with $\codim_X\Sing\,X\gess 5$.
\sk
In Section~1 we review some basics of filtered $\D$-modules, filtered differential complexes, and Du~Bois complexes.
In Section~2 we prove Theorem~\hl{T1}{1} after reviewing Du~Bois complexes in the hypersurface case.
In Section~3 we prove Theorems~\hl{T2}{2}--\hl{T3}{3} and Corollary~\hl{C3}{3} using the nearby and vanishing cycle Hodge modules.
\msn
{\bf Convention.} In this paper a complex algebraic variety means a separated reduced scheme of finite type over $\C$, although we consider only {\it closed points\1} (except in 2.3). So it is a variety in the sense of Serre \cite[\S 34]{Se} (that is, {\it irreducibility\1} is not supposed). Sometimes $\Q_{Y^{\rm an}}$ is denoted by $\Q_Y$ to simplify the notation.
\msn
{\bf Acknowledgements.} We thank the referee for useful comments to improve the paper.
\bs\bs
\vbox{\centerline{\bf 1. Filtered differential complexes and Du~Bois complexes}
\bsn
In this section we review some basics of filtered $\D$-modules, filtered differential complexes, and Du~Bois complexes.}
\msn
{\bf 1.1.~Filtered $\D$-modules and filtered differential complexes.} Let $Y$ be a smooth complex algebraic variety. Let $\DbcF(\D_Y)$ be the derived category of bounded complexes of filtered left $\D_Y$-modules with coherent cohomology sheaves. Here the cohomology sheaves are taken in the category of graded $\mopl_p\,F_p\D_Y$-modules, see \cite{mhp}.
We denote by $\DbcF(\D_Y)^r$ a similar category for filtered {\it right\1} $\D_Y$-modules. There is an equivalence of categories
$$\begin{array}{rl}\tau_Y:\DbcF(\D_Y)&\!\!\!\!\simto\DbcF(\D_Y)^r\\ {\rotatebox{90}{$\in$}}\q\q&\q\q\q\q{\rotatebox{90}{$\in$}}\\ (M^{\ssb},F)\,\,\,&\!\!\!\mapsto\,\,\,(\Om^{d_Y}_Y,F){\otimes_{\OO_Y}}(M^{\ssb},F),\end{array}$$
\vskip-5mm
\htt{1.1.1}{}
$$\h{with}\q\q\q\Gr^F_p\Om^{d_Y}_Y\eq0\q\q(p\ne-d_Y)\raise7mm\h{}.
\leqno(1.1.1)$$
Note that the filtration $F$ on $\Om^{d_Y}_Y$ is {\it shifted by\1} $-d_Y$. In this paper we {\it distinguish\1} $\Om^{d_Y}_Y$ with the dualizing sheaf $\om_Y$ on which the filtration $F$ is {\it not\1} shifted, see (\hl{1.4.1}{1.4.1}) below and also \cite[(B.2.2--3)]{ypg}.
\sk
Let $\DbcF(\OO_Y,\Diff)$ be the derived category of bounded filtered differential complexes $(L^{\ssb},F)$. Here each $(L^j,F)$ is a filtered $\OO_Y$-module, and the differential consists of filtered differential morphisms. We assume that the $\Hc^j\Gr_F^pL^{\ssb}$ are coherent $\OO_Y$-modules for any $p,j\in\Z$, and $F$ is a bounded filtration. We usually use {\it decreasing\1} filtration $F^{\ssb}$ for a filtered differential complexes. However, {\it increasing\1} filtration $F_{\ssb}$ is better for the relation with filtered $\D$-modules. These are related by $F_p=F^{-p}$ in general.
\sk
Recall that $\phi:(L^{\ssb},F)\to(L'{}^{\ssb},F)$ is a {\it filtered differential morphism\1} (in the sense of \cite[2.2.4]{mhp}) if the composition
$$F_pL\into L\buildrel{\phi}\over\to L'\onto L'/F_{p-q-1}L'$$
is a differential morphism of order $\les q$ for any $p,q\in\Z$. (Note that a differential morphism of strictly negative order is 0.) We say that $\phi$ is a {\it differential morphism\1} if it belongs to the image of the {\it injective\1} morphism
\htt{1.1.4}{}
$${\rm Hom}_{\D_Y}(L{\otimes}_{\OO_Y}\D_Y,L'{\otimes}_{\OO_Y}\D_Y)\into{\rm Hom}_{\C}(L,L'),
\leqno(1.1.4)$$
induced by ${\otimes}_{\D_Y}\OO_Y$. Here the source has the filtration $F$ induced from the one on $\D_Y$ via the isomorphism
$${\rm Hom}_{\D_Y}(L{\otimes}_{\OO_Y}\D_Y,L'{\otimes}_{\OO_Y}\D_Y)={\rm Hom}_{\OO_Y}(L,L'{\otimes}_{\OO_Y}\D_Y).$$
We say that a differential morphism has {\it order} $\les p$ if it belongs to the image of $F_p$ of the left-hand side of the above isomorphism, see \cite{mhp}, \cite{ind}.
\msn
{\bf 1.2.~Equivalence of categories.} There are functors
\htt{1.2.1}{}
$$\aligned\DR_Y^{-1}&:\DbcF(\OO_Y,\Diff)\to\DbcF(\D_Y)^r,\\
\DR_Y&:\DbcF(\D_Y)\to\DbcF(\OO_Y,\Diff),\\
\DR_Y^r&:\DbcF(\D_Y)^r\to\DbcF(\OO_Y,\Diff),\endaligned
\leqno(1.2.1)$$
inducing {\it equivalences of categories.} Indeed, there are canonical isomorphisms in the derived categories
\htt{1.2.2}{}
$$\aligned&\DR_Y^r\ssc\tau_Y=\DR_Y,\\&\DR_Y^{-1}\ssc\DR_Y^r\cong{\rm id},\\&\DR_Y^r\ssc\DR_Y^{-1}\cong{\rm id}.\endaligned
\leqno(1.2.2)$$
Here $\DR_Y$ is the usual de Rham functor (shifted by $d_Y$ to the left as usual). The $i$\1th component of $\DR_Y(M,F)$ for a filtered left $\D_Y$-module $(M,F)$ is given by
\htt{1.2.3}{}
$$\Om_Y^{i+d_Y}{\otimes}_{\OO_Y}(M,F[-i{-}d_Y])\q(i\in[-d_Y,0]).
\leqno(1.2.3)$$
The functor $\DR_Y^{-1}$ is defined by ${\otimes}_{\OO_X}(\D_X,F)$ using the {\it injectivity\1} of (\hl{1.1.4}{1.1.4}).
For $\DR_Y^r$, we need the sheaf of vector fields $\Theta_Y$. The $i$\1th component of $\DR_Y^r({}^r\!M,F)$ for a filtered right $\D_Y$-module $({}^r\!M,F)$ is given by
\htt{1.2.4}{}
$$({}^r\!M,F[-i]){\otimes}_{\OO_Y}\h{$\bigwedge$}^{-i}\Theta_Y\q(i\in[-d_Y,0]),
\leqno(1.2.4)$$
where $(F[m])_p=F_{p-m}$ in a compatible way with $(F[m])^p=F^{p+m}$ and $F^p=F_{-p}$.
\sk
The first isomorphism in (\hl{1.2.2}{1.2.2}) is an isomorphism of filtered complexes, and the second and third ones are canonical filtered quasi-isomorphisms. The proofs are reduced to the calculation of the Koszul complex for the actions of $\xi_1,\dots,\xi_{d_Y}$ on $\OO_Y[\xi_1,\dots,\xi_{d_Y}]$, where $\xi_i:=\Gr^F_1\dd_{y_i}$ with $y_1,\dots,y_{d_Y}$ local coordinates of $Y$, see \cite{mhp}, \cite{ind}. (Concerning the coherence, note that $\OO_{Y,y}[\xi_1,\dots,\xi_{d_Y}]$ is a regular ring, and any finite graded module over it has a free resolution of finite length.)
\msn
{\bf 1.3.~Direct image functor.} Let $f:X\to Y$ be a {\it proper\1} morphism of smooth complex algebraic varieties. We have the direct image (or pushforward) functors
\htt{1.3.1}{}
$$\aligned f_*^\Diff&:\DbcF(\OO_X,\Diff)\to\DbcF(\OO_Y,\Diff),\\ f_*^{\D}&:\DbcF(\D_X)\to\DbcF(\D_Y),\\ f_*^{\D}&:\DbcF(\D_X)^r\to\DbcF(\D_Y)^r.\endaligned
\leqno(1.3.1)$$
The direct image functor $f_*^\Diff$ is defined as the {\it sheaf-theoretic direct image\1} taking a canonical flasque resolution \cite{Go} (together with the filtered truncation $\tau_{\les k}$ for $k\gg 0$), since the differential morphisms are stable by direct images. The remaining functors $f_*^{\D}$ are defined so that they are compatible with $f_*^\Diff$ via the equivalence of categories induced by the functors in (\hl{1.2.1}{1.2.1}). It is easy to see that the functor $f_*^{\D}$ is canonically isomorphic to the usual direct image of filtered left or right $\D$-modules, see also \cite{ypg}.
\sk
These direct image functors are compatible with the direct image functor of mixed Hodge modules (forgetting the weight filtration) only in the case the morphism $f$ is {\it proper.}
\msn
{\bf 1.4.~Dual functor.} Let $\om_Y\to\I_Y^{\ssb}$ be an injective resolution as a right $\D_Y$-module. There are commuting two structures of right $\D_Y$-modules on $\om_Y{\otimes}_{\OO_Y}\D_Y$, and on $\I_Y^j{\otimes}_{\OO_Y}\D_Y$. We have the filtration $F$ on $\om_Y$, $\I_Y^j$ defined by
\htt{1.4.1}{}
$$\Gr_F^p\om_Y\eq\Gr_F^p\I_Y^j\eq0\q(p\ne 0).
\leqno(1.4.1)$$
We then get a filtered quasi-isomorphism of complex of modules having commuting two right $\D_Y$-module structures (using $F_p=F^{-p}$)
$$(\om_Y^{\ssb},F){\otimes}_{\OO_Y}(\D_Y,F)\simto(\I_Y^{\ssb},F){\otimes}_{\OO_Y}(\D_Y,F).$$
\sk
We can define the dual $\DD(L^{\ssb},F)$ for $(L^{\ssb},F)\in\DbcF(\OO_Y,\Diff)$ by
\htt{1.4.2}{}
$$\DD(L,F):=\tau_{\les k}\Hc om_{\OO_Y}\bl((L^{\ssb},F),(\I_Y^{\ssb}[d_Y],F)\br),
\leqno(1.4.2)$$
for $k\gg 0$. The differential of complex is given so that we have the canonical isomorphism of filtered complexes
\htt{1.4.3}{}
$$\aligned&\DR_Y^{-1}\bl(\Hc om_{\OO_Y}\bl((L^{\ssb},F),(\I_Y^{\ssb}[d_Y],F)\br)\br)\\={}&\Hc om_{\D_Y}\bl(\DR_Y^{-1}(L^{\ssb},F),(\I_Y^{\ssb}[d_Y],F){\otimes}_{\OO_Y}(\D_Y,F)\br),\endaligned
\leqno(1.4.3)$$
where the right $\D_Y$-module structure associated with $\I_Y^{\ssb}{\otimes}_{\OO_Y}$ is used for $\Hc om_{\D_Y}$.
Using the equivalence of categories induced by (\hl{1.2.1}{1.2.1}), this construction implies the dual functor $\DD$ also for $\DbcF(\D_Y)^r$, $\DbcF(\D_Y)$ together with the canonical isomorphisms
\htt{1.4.4}{}
$$\aligned\DR_Y^{-1}\ssc\DD&=\DD\ssc\DR_Y^{-1},\\ \DR_Y^r\ssc\DD&=\DD\ssc\DR_Y^r,\\ \tau_Y\ssc\DD&=\DD\ssc\tau_Y.\endaligned
\leqno(1.4.4)$$
By (\hl{1.4.1}{1.4.1}) we have also the canonical isomorphism
\htt{1.4.5}{}
$$\Gr_F^p\DD(L^{\ssb},F)=\DD(\Gr_F^{-p}L^{\ssb}),
\leqno(1.4.5)$$
since the $\I_Y^j$ are injective $\OO_Y$-modules (using the flatness of $\D_Y$ over $\OO_Y$).
Note that for $N^{\ssb}\in D^b_{\rm coh}(\OO_Y)$ in general, we have by definition
\htt{1.4.6}{}
$$\DD(N^{\ssb}):=\RR\Hc om_{\OO_Y}(N^{\ssb},\om_Y[d_Y]).
\leqno(1.4.6)$$
\sk
These functors $\DD$ are compatible with the functor $\DD$ for mixed Hodge modules (forgetting the weight filtration).
\msn
{\bf 1.5.~Du~Bois complexes.} Let $X$ be a reduced complex algebraic variety (where $X$ is not assumed irreducible). We have the filtered Du~Bois complex $(\Omu_X^{\ssb},F)\in\DbcF(\OO_Y,\Diff)$ using a canonical functor from the derived category of filtered bounded differential complexes with order of differential at most 1 in \cite{DB} to $D^bF(\OO_Y,\Diff)$, see also \cite{Fi}. Using a simplicial resolution $\gamma:Y_{\ssb}\to X$, it can be defined as
\htt{1.5.1}{}
$$(\Omu^{\ssb},F):=\RR\gamma_*(\Om_{Y_{\bullet}}^{\ssb},F).
\leqno(1.5.1)$$
This is a single complex associated with a triple complex whose $(i,j,k)$ component is $\gamma_*\J^k(\Om_{Y_i}^j)$, where $\J^k(M)$ denotes the $k$\1th component of the canonical flasque resolution of a sheaf $M$ in general, see \cite{Go}. The Hodge filtration is induced by the so-called ``stupid truncation" $\sigma_{\ges p}$ on the $\Om_{Y_i}^{\ssb}$. (For a complex $K^{\ssb}$ in general, we have $(\sigma_{\ges p}K)^j=K^j$ if $j\gess p$, and 0 otherwise.) We then get the isomorphism
\htt{1.5.2}{}
$$\Omu_X^p=\RR\gamma_*\Om_{Y_{\bullet}}^p,
\leqno(1.5.2)$$
where $\Omu_X^p:=\Gr_F^p\Omu_X^{\ssb}[p]\in D^b_{\rm coh}(\OO_X)$. Note that $\Omu_X^p$ does {\it not\1} mean the $p$\1th component of the complex $\Omu_X^{\ssb}$; the latter does not make sense for an object of a derived category, and $\Omu_X^p$, which is an abbreviation of $(\Omu_X^p)^{\ssb}$, is still a complex on $X$ in general.
\sk
We can also use a hypercubic resolution as in \cite{GNPP} (instead of a simplicial resolution). This gives an isomorphism
\htt{1.5.3}{}
$$\Omu_X^{d_X}=\mopl_i\,(\gamma_i)_*\Om_{\Xt_i}^{d_X},
\leqno(1.5.3)$$
using the Grauert-Riemenschneider vanishing theorem \cite{GR} for the higher direct images. Here the $\gamma_i:\Xt_i\to X_i$ are desingularizations of the irreducible components $X_i$ of $X$ with $d_{X_i}\eq d_X$, see also \cite{kol}.
\sk
As an immediate corollary of \cite[Thm.~4.2]{mhc}, we have the following.
\par\htt{T1.5}{}\msn
{\bf Theorem~1.5.} {\it Assume $X$ is a closed subvariety of a smooth complex algebraic variety $Y$ $($for instance, $X$ is quasi-projective$)$ so that there is $(M_X^{\prime\ssb},F)\in\DbcF(\D_Y)$ underlying $\Q_{h,X}\in D^b{\rm MHM}(Y)$. Then we have a canonical isomorphism}
\htt{1.5.4}{}
$$\DR_Y(M_X^{\prime\ssb},F)=(\Omu_X^{\ssb},F)\q\h{in}\,\,\,\,\DbcF(\OO_Y,\Diff).
\leqno(1.5.4)$$
\par\htt{R1.5a}{}\msn
{\bf Remark~1.5a.} If $X$ is embeddable into a smooth variety $Y$, it is rather easy to show a {\it non-canonical\1} isomorphism in (\hl{1.5.4}{1.5.4}) as follows (compare to \cite{MOPW}). Taking an embedded resolution $\pi:(\Yt,E)\to(Y,X)$,
there is a distinguished triangle in $\DbcF(\OO_Y,\Diff)$\,:
\htt{1.5.5}{}
$$(\Omu_Y^{\ssb},F)\to(\Omu_X^{\ssb},F)\oplus\RR\pi^{\Diff}_*(\Omu_{\Yt}^{\ssb},F)\to\RR\pi^{\Diff}_*(\Omu_E^{\ssb},F)\buildrel{+1}\over\to,
\leqno(1.5.5)$$
see \cite[Prop.~4.11]{DB}. This is equivalent to the filtered acyclicity of the single complex associated with the commutative diagram in $\DbcF(\OO_Y,\Diff)$\,:
\htt{1.5.6}{}
$$\begin{array}{ccccc}\RR\pi^{\Diff}_*(\Omu_{\Yt}^{\ssb},F)&\to&\RR\pi^{\Diff}_*(\Omu_E^{\ssb},F)\\ \uparrow&\raise4mm\h{}\raise-3mm\h{}&\uparrow\\(\Omu_Y^{\ssb},F)&\to&(\Omu_X^{\ssb},F)\end{array}
\leqno(1.5.6)$$
(This is not necessarily a triple complex, since the diagram commutes up to a {\it homotopy.}
We must use this homotopy in the definition of the differential of the associated single complex. Here we can argue in $\DbcF(\D_Y)$ by applying the functor $\DR_Y^{-1}$.) We then see that (\hl{1.5.5}{1.5.5}) and (\hl{1.5.6}{1.5.6}) are equivalent to the isomorphism in $\DbcF(\OO_Y,\Diff)$\,:
\htt{1.5.7}{}
$$\aligned&C\bl((\Omu_Y^{\ssb},F)\to(\Omu_X^{\ssb},F)\br)\\\simto{}&C\bl(\RR\pi^{\Diff}_*(\Omu_{\Yt}^{\ssb},F)\to\RR\pi^{\Diff}_*(\Omu_E^{\ssb},F)\br).\endaligned
\leqno(1.5.7)$$
(This is {\it not\1} canonical, since it depends on the choice of the homotopy explained above.)
\sk
Let $\bl(M^{\ssb}_Y(!X),F\br)\in\DbcF(\D_Y)$ be the underlying complex of filtered left $\D_Y$-modules of $(j_U)_*\Q_{h,U}\in D^b{\rm MHM}(Y)$ with $j_U:U:=Y\,{\setminus}\,X\into Y$ the inclusion morphism. Set
$$\bl(\Omu_Y^{\ssb}(!X),F\br):=\DR_Y\bl(M^{\ssb}_Y(!X),F\br)\in\DbcF(\OO_Y,\Diff).$$
Define similarly $\bl(M^{\ssb}_{\Yt}(!E),F\br)$, $\bl(\Omu_{\Yt}^{\ssb}(!E),F\br)$ replacing $(Y,X)$ with $(\Yt,E)$. Since $\pi$ is proper and $\pi\ssc\jt_U\eq j_U$ with $\jt_U:U\into\Yt$ the natural inclusion, we get the canonical isomorphisms
\htt{1.5.8}{}
$$\aligned\RR\pi_*^{\D}\bl(M^{\ssb}_{\Yt}(!E),F\br)&=\bl(M^{\ssb}_Y(!X),F\br)\q\h{in}\,\,\,\DbcF(\D_Y),\\ \RR\pi_*^{\Diff}\bl(\Omu_{\Yt}^{\ssb}(!E),F\br)&=\bl(\Omu_Y^{\ssb}(!X),F\br)\q\h{in}\,\,\,\DbcF(\OO_Y,\Diff).\endaligned
\leqno(1.5.8)$$
Since $E$ is a divisor with normal crossings on $\Yt$, it is well-known that $\DR_{\Yt}^{-1}(\Omu_E^{\ssb},F)[d_E]$ is isomorphic to the underlying filtered $\D_Y$-module of the mixed Hodge module $\Q_{h,E}[d_E]$. Combining this with (\hl{1.5.8}{1.5.8}), we see that the right-hand side of (\hl{1.5.7}{1.5.7}) is isomorphic to $\bl(\Omu_Y^{\ssb}(!X),F\br)[1]$. So we get the non-canonical isomorphism
\htt{1.5.9}{}
$$C\bl((\Omu_Y^{\ssb},F)\to(\Omu_X^{\ssb},F)\br)\cong\bl(\Omu_Y^{\ssb}(!X),F\br)[1],
\leqno(1.5.9)$$
which is equivalent to a non-canonical isomorphism
\htt{1.5.10}{}
$$C\bl(\bl(\Omu_Y^{\ssb}(!X),F\br)\to(\Omu_Y^{\ssb},F)\br)\cong(\Omu_X^{\ssb},F).
\leqno(1.5.10)$$
Here the morphism $\bl(\Omu_Y^{\ssb}(!X),F\br)\to(\Omu_Y^{\ssb},F)$, or equivalently, $M^{\ssb}_Y(!X)\to(\OO_Y,F)[-d_Y]$, is {\it unique.} Indeed, the filtration $F$ on $M^{\ssb}_Y(!X)$ is {\it strict,} and we can apply the vanishing of negative extension groups by using the canonical truncations $\tau_{\les k}$ together with graded $\mopl_p\,F_p\D_Y$-modules. (Note that $\Hc^jM^{\ssb}_Y(!X)\eq 0$, that is, $H^j(j_U)_!\Q_{h,U}\eq 0$, if $j>d_Y$.) The non-canonical isomorphism (\hl{1.5.4}{1.5.4}) then follows.
\par\htt{R1.5b}{}\msn
{\bf Remark~1.5b.} In the notation and assumption of Remark~\hl{R1.5a}{1.5a} above, assume $\Q_{h,X}[d_X]$ is a mixed Hodge module, or equivalently, $M_X:=M_X^{\prime\ssb}[d_X]$ is a $\D_Y$-module; for instance, $X$ is a hypersurface or more generally a local complete intersection. Then the isomorphism (\hl{1.5.4}{1.5.4}) is {\it unique\1} using the injectivity of
\htt{1.5.11}{}
$${\rm End}_{\D_Y}(M_X,F)\into{\rm End}_{\D_Y}(M_X).
\leqno(1.5.11)$$
Indeed, we have the isomorphisms
\htt{1.5.12}{}
$${\rm End}_{\D_Y}(M_X)={\rm End}_{\C_Y}(\C_X[d_X])=\C,
\leqno(1.5.12)$$
in the $X$ connected case, since $M_X$ is {\it regular holonomic.} We can then restrict to the smooth part of $X$.
\msn
{\bf 1.6.~K\"ahler differential forms.} Let $X$ be a reduced closed subvariety of a smooth complex algebraic variety $Y$. Let $\I_X\subset\OO_Y$ be the reduced ideal sheaf of $X$. The sheaf $\Om_X^p$ of K\"ahler differential forms of degree $p$ can be defined as the quotient $\OO_X$-module of $\Om_Y^p|_X$ divided by the $\OO_X$-submodule generated locally by $\ddd g\,{\wedge}\,\Om_Y^{p-1}|_X$ for $g\in\I_X$, see for instance \cite[Def.\,1.2]{GK}, \cite[Def.\,1]{Re}. (Here the restriction $|_X$ means the tensor product with $\OO_X$ over $\OO_Y$.)
\sk
If $X$ is a reduced hypersurface defined by a function $f$, we have
\htt{1.6.1}{}
$$\Om_X^p=(\Om_Y^p/\ddd f{\wedge}\,\Om_Y^{p-1})|_X\q\q(p\in\Z).
\leqno(1.6.1)$$
\par\htt{R1.6}{}\msn
{\bf Remark~1.6.} It is well-known that $\Om_X^p$ may have {\it torsion.} In the case $Y\eq\C^2$ with $f\eq y_1y_2$, for instance, it is easy to see that
\htt{1.6.2}{}
$$\Gamma(X,\Om_X^1)={\rm Coker}\bl(y_1y_2:(y_1,y_2)\to(y_1,y_2)\br),
\leqno(1.6.2)$$
with $(y_1,y_2)\subset\C[y_1,y_2]$ the maximal ideal. Indeed, the partial derivatives $\dd_{y_i}f$ form a regular sequence so that
\htt{1.6.3}{}
$$\aligned\Gamma(Y,\Om_Y^1)/\ddd f{\wedge}\1\Gamma(Y,\OO_Y)&=\ddd f{\wedge}\,\Gamma(Y,\Om_Y^1)\\&\cong(y_1,y_2)\subset\C[y_1,y_2].\endaligned
\leqno(1.6.3)$$
\sk
This can be extended to the case $f\eq y_1y_2\in\C[y_1,\dots,y_n]$ with $Y\eq\C^n$ ($n\gess 3$). Here we have the decompositions $Y\eq Y_1\,{\times}\,Y_2$, $X\eq X_1\,{\times}\,X_2$ with $Y_1\eq\C^2$, $X_2\eq Y_2\eq\C^{n-2}$, and
\htt{1.6.4}{}
$$\Om_X^1={\rm pr}_1^*\Om_{X_1}^1\oplus{\rm pr}_2^*\Om_{X_2}^1,
\leqno(1.6.4)$$
with ${\rm pr}_i:X\to X_i$ the canonical projections.
\bs\bs
\vbox{\centerline{\bf 2. Proof of Theorem~\hl{T1}{1}}
\bsn
In this section we prove Theorem~\hl{T1}{1} after reviewing Du~Bois complexes in the hypersurface case.}
\msn
{\bf 2.1.~Du~Bois complexes of hypersurfaces.} Assume $X$ is a reduced hypersurface of a smooth complex algebraic variety $Y$ with $i_X:X\into Y$, $j_U:U:=Y\,{\setminus}\,X\into Y$ the inclusion morphisms. We have a filtered left $\D_Y$-module $(\OO_Y(*X),F)$ underlying the mixed Hodge module $(j_U)_*\Q_{h,U}[d_Y]$. This corresponds to the filtered right module $(\Om^{d_Y}_Y(*X),F)$, where the filtration is {\it shifted by\1} $-d_Y$ so that $F_p\1\Om_Y^{d_Y}(*X)\ne 0\iff p\ges-d_Y$, see also a remark after (\hl{1.1.1}{1.1.1}).
\sk
Consider the quotient filtered right $\D_Y$-module $\bl(\Om^{d_Y}_Y(*X)/\Om^{d_Y}_Y,F\br)$. This underlies a mixed Hodge module
$$H^1i_X^!(\Q_{h,Y}[d_Y])=i_X^!\Q_{h,Y}[d_Y{+}1],$$
using the distinguished triangle
\htt{2.1.1}{}
$$(i_X)_*i_X^!\Q_{h,Y}\to\Q_{h,Y}\to(j_U)_*\Q_{h,U}\buildrel{+1}\over\to.
\leqno(2.1.1)$$
\sk
We have a canonical isomorphism of mixed Hodge modules
\htt{2.1.2}{}
$$\DD\1i_X^!(\Q_{h,Y}[d_Y{+}1])=\Q_{h,X}(d_Y)[d_X],
\leqno(2.1.2)$$
since
\htt{2.1.3}{}
$$\DD\ssc i_X^!\eq i_X^*\ssc\DD,\q\DD(\Q_{h,Y}[d_Y])=\Q_{h,Y}(d_Y)[d_Y].
\leqno(2.1.3)$$
Passing to the underlying filtered $\D$-modules, this implies the canonical isomorphism of filtered right $\D_Y$-modules
\htt{2.1.4}{}
$$\DD\bl(\Om^{d_Y}_Y(*X)/\Om^{d_Y}_Y,F\br)=({}^r\!M_X,F[d_Y]),
\leqno(2.1.4)$$
where $({}^r\!M_X,F):=\tau_Y(M_X,F)$ with $(M_X,F)$ the underlying filtered left $\D_Y$-module of $\Q_{h,X}[d_X]$. Set
\htt{2.1.5}{}
$$\bl(\om_Y(*X),F\br):=\bl(\Om^{d_Y}_Y(*X),F[d_Y]\br),
\leqno(2.1.5)$$
This is compatible with $(\om_Y,F)\eq(\Om^{d_Y}_Y,F[d_Y])$, see (\hl{1.1.1}{1.1.1}), (\hl{1.4.1}{1.4.1}), and \cite[(B.2.2--3)]{ypg}. The isomorphism (\hl{2.1.4}{2.1.4}) is then equivalent to the canonical isomorphism
\htt{2.1.6}{}
$$\DD\bl(\om_Y(*X)/\om_Y,F\br)=({}^r\!M_X,F).
\leqno(2.1.6)$$
\sk
Note that the shift of filtration $[m]$ for $m\in\Z$ changes the weight like the {\it Tate twist\1} $(m)$, that is, {\it the weight is decreased by\1} $2m$. So $(\om_Y,F)$ and $\bl(\om_Y(*X)/\om_Y,F\br)$ are pure of weight $-d_Y$ and $1{-}d_Y\eq{-}d_X$ respectively.
\sk
The functor $\DD$ is compatible with the de Rham functor $\DR_Y^r$, see (\hl{1.4.4}{1.4.4}). So (\hl{2.1.6}{2.1.6}) together with (\hl{1.4.5}{1.4.5}) and Theorem~\hl{T1.5}{1.5} implies the following (compare to \cite{MOPW}).
\par\htt{P2.1}{}\msn
{\bf Proposition~2.1.} {\it There are canonical isomorphisms in $D^b_{\rm coh}(\OO_Y)\,{:}$
\htt{2.1.7}{}
$$\aligned\DD\Gr^F_p\DR_Y^r\bl(\om_Y(*X)/\om_Y\br)&=\Gr_F^p\DR_Y(M_X)\\&=\Omu_X^p[d_X{-}p]\q\q(p\in\Z),\endaligned
\leqno(2.1.7)$$
where $F_p\eq F^{-p}$.}
\par\htt{R2.1}{}\msn
{\bf Remark~2.1.} The canonical isomorphism (\hl{2}{2}) in $D^b_{\rm coh}(\OO_X)$ is equivalent to a canonical isomorphism of $\OO_X$-modules
\htt{2.1.8}{}
$$\Om_X^p\simto\Hc^0(\Omu_X^p)\q\h{together with}\q\Hc^j(\Omu_X^p)=0\,\,\,\,(j\ne 0).
\leqno(2.1.8)$$
Here the first morphism is induced from a property of K\"ahler differentials.
\msn
{\bf 2.2.~Proof of the implication} {\rm (\hl{c}{c}) $\Rightarrow$ (\hl{b}{b})}. Assume condition~(\hl{c}{c}). We have the isomorphisms for $q\in[0,p]$\,:
\htt{2.2.1}{}
$$\Gr^F_q\DR_Y^r(\om_Y(*X)/\om_Y)=\bl(\sigma_{\ges d_Y-q}(\Om_Y^{\ssb}/f\Om_Y^{\ssb},\ddd f\wedge)\br)[d_Y],
\leqno(2.2.1)$$
using (\hl{2.1.5}{2.1.5}) and (\hl{5}{5}) in the introduction. Here $\sigma_{\ges k}$ is the ``stupid truncation", see a remark before (\hl{1.5.2}{1.5.2}). Indeed, the differential of $\Om_Y^{\ssb}(*X)$ is given by
\htt{2.2.2}{}
$$\ddd(f^{-q-1}\eta)=-(q{+}1)\1f^{-q-2}\ddd f{\wedge}\eta+f^{-q-1}\ddd\eta\q(\eta\in\Om_Y^k),
\leqno(2.2.2)$$
where the last term can be neglected by passing to the graded pieces (since (\hl{c}{c}) is assumed). The multiplication by $-(q{+}1)$ can be adjusted using appropriate scalar multiplications on the components of complex. We may assume that this scalar multiplication is the identity on the component of degree $-q$.
\sk
The right-hand side of (\hl{2.2.1}{2.2.1}) is quasi-isomorphic to the shifted mapping cone
\htt{2.2.3}{}
$$C\bl(\sigma_{\ges d_Y-q}(\Om_Y^{\ssb},\ddd f\wedge)\buildrel{f\,}\over\to\sigma_{\ges d_Y-q}(\Om_Y^{\ssb},\ddd f\wedge)\br)[d_Y],
\leqno(2.2.3)$$
which is a complex of locally free $\OO_Y$-modules.
\sk
Applying the functor $\DD$ as in (\hl{1.4.6}{1.4.6}), and using the canonical isomorphisms
$$\Hc om_{\OO_Y}(\Om_Y^k,\om_Y)=\Om_Y^{d_Y-k}\q\q(k\in\Z),$$
we then get by Proposition~\hl{P2.1}{2.1} the isomorphisms for $q\in[0,p]$\,:
\htt{2.2.4}{}
$$\Omu_X^q[d_X{-}q]=C\bl(\sigma_{\les q}(\Om_Y^{\ssb},\ddd f\wedge)\buildrel{f\,}\over\to\sigma_{\les q}(\Om_Y^{\ssb},\ddd f\wedge)\br)[d_Y{-}1].
\leqno(2.2.4)$$
\sk
By Proposition~\hl{P1}{1} together with \cite[Cor.\,1]{hi} for (\hl{a}{a})\,$\Leftrightarrow$\,(\hl{c}{c}), we get the strict inequality
\htt{2.2.5}{}
$$\codim_Y\Sing\,X>p{+}2,
\leqno(2.2.5)$$
since these imply that $\tfrac{1}{2}\1\codim_Y\Sing\,X\ges\alt_X\ges p{+}1\ges 2$.
\sk
Combined with Proposition~\hl{P2}{2}, this gives the quasi-isomorphisms for $q\in[0,p]$\,:
\htt{2.2.6}{}
$$\sigma_{\les q}(\Om_Y^{\ssb},\ddd f\wedge)\simto\Om_Y^q/\ddd f{\wedge}\,\Om_Y^{q-1}[-q].
\leqno(2.2.6)$$
\sk
In order to get the isomorphisms
\htt{2.2.7}{}
$$\Om_X^q\cong\Omu_X^q\q(q\in[0,p]),
\leqno(2.2.7)$$
it is then enough to show that $\Om_Y^q/\ddd f{\wedge}\Om_Y^{q-1}$ is $f$-torsion-free in view of the definition of $\Om_X^q$, see (\hl{1.6.1}{1.6.1}). But Proposition~\hl{P2}{2} and (\hl{2.2.5}{2.2.5}) imply the isomorphisms for $q\in[0,p]$\,:
\htt{2.2.8}{}
$$\Om_Y^q/\ddd f{\wedge}\,\Om_Y^{q-1}=\ddd f{\wedge}\,\Om_Y^q\,\,\bl(\subset\Om_Y^{q+1}\br),
\leqno(2.2.8)$$
and the right-hand side is $f$-torsion-free. So the isomorphisms in (\hl{2.2.7}{2.2.7}) follow.
\sk
It remains to show that the obtained isomorphism (\hl{2.2.7}{2.2.7}) coincides with the {\it canonical\1} morphism. We first prove this in the $X$ {\it smooth\1} case. There are canonical isomorphisms in this case for any $q\in\Z$\,:
\htt{2.2.9}{}
$$\aligned\Gr&^F_q\DR_Y^r\bl(\om_Y(*X)/\om_Y\br)[-q]\\&{=}\,\,{\rm Ker}\bl(-\ddd f/f\,\wedge:\Om_Y^{d_Y-q}(X)|_X\to\Om_Y^{d_Y-q+1}(2X)|_X\br)\\(&{=}\,\,\Om_X^{d_X-q}).\endaligned
\leqno(2.2.9)$$
Here $|_X$ means the tensor product with $\OO_X$ over $\OO_Y$, and a remark after (\hl{2.2.2}{2.2.2}) is used. The second term of (\hl{2.2.9}{2.2.9}) is independent of the choice of $f$. (The last isomorphism can be obtained by taking the {\it residue.})
\sk
Consider the canonical short exact sequences
\htt{2.2.10}{}
$$0\to\Om_Y^{d_Y-q}\to\Om_Y^{d_Y-q}(X)\to\Om_Y^{d_Y-q}(X)|_X\to0.
\leqno(2.2.10)$$
Applying the functor $\DD$, we get short exact sequences which are naturally isomorphic to
\htt{2.2.11}{}
$$0\to\Om_Y^q(-X)\to\Om_Y^q\to\Om_Y^q|_X\to0.
\leqno(2.2.11)$$
The above construction of (\hl{2.2.7}{2.2.7}) is then identified with the canonical morphism. So the proof of the implication (\hl{c}{c}) $\Rightarrow$ (\hl{b}{b}) is reduced to the following (here we consider the image of the difference of two morphisms, whose support is contained in the singular locus of $X$).
\par\htt{P2.2}{}\msn
{\bf Proposition~2.2.} {\it Assume a hypersurface $X\subset Y$ is higher $p$-log canonical for some positive integer $p$. Then, for $q\in[0,p]$, the coherent sheaves of K\"ahler differential forms $\Om_X^q$ are torsion-free\,$;$ more precisely, we have
\htt{2.2.12}{}
$$\Hc^0_{\Sing\,X}\Om_X^q\eq 0\q(q\in[0,p]),
\leqno(2.2.12)$$
so the restriction morphisms ${\rm End}_{\OO_X}(\Om_X^q)\,{\to}\,{\rm End}_{\OO_{X'}}(\Om_{X'}^q)$ with $X'\,{:=}\,X\1{\setminus}\1\Sing\1 X$ are injective.}
\msn
{\it Proof.} We can get the inclusion
\htt{2.2.13}{}
$$\Om_X^q\into\Om_Y^{q+1}|_X,
\leqno(2.2.13)$$
applying the snake lemma to the diagram
\htt{2.2.14}{}
$$\begin{array}{ccccccccc}0&\to&\Om_Y^q/\ddd f{\wedge}\Om_Y^{q-1}&\to&\Om_Y^{q+1}&\to&\Om_Y^{q+1}/\ddd f{\wedge}\Om_Y^q&\to&0\\&&\downarrow\!\raise.5mm\h{$\scriptstyle f$}&&\downarrow\!\raise.5mm\h{$\scriptstyle f$}&\raise4mm\h{}\raise-2.3mm\h{}&\downarrow\!\raise.5mm\h{$\scriptstyle f$}\\0&\to&\Om_Y^q/\ddd f{\wedge}\Om_Y^{q-1}&\to&\Om_Y^{q+1}&\to&\Om_Y^{q+1}/\ddd f{\wedge}\Om_Y^q&\to&0\end{array}
\leqno(2.2.14)$$
and using Proposition~\hl{P2}{2} and (\hl{2.2.5}{2.2.5}). Here $\Om_Y^{q+1}/\ddd f{\wedge}\Om_Y^q$ is $f$-torsion-free by an isomorphism as in (\hl{2.2.8}{2.2.8}), since $q{+}1<p{+}2<\codim_Y\,\Sing\,X$.
Proposition~\hl{P2.2}{2.2} thus follows.
\par\htt{R2.2}{}\msn
{\bf Remark~2.2.} In the definition of higher $p$-Du~Bois singularity for hypersurfaces, it is enough to assume a {\it non-canonical\1} isomorphism (\hl{2}{2}) locally on $X$. This follows from the proof of the implication (\hl{b}{b}) $\Rightarrow$ (\hl{c}{c}) just below, where it is not needed that the isomorphism (\hl{2}{2}) is a canonical one.
\msn
{\bf 2.3.~Proof of the implication} {\rm (\hl{b}{b}) $\Rightarrow$ (\hl{c}{c})}. We argue by increasing induction on $p\ges 1$. Assume condition~(\hl{b}{b}) for $p$ holds, but (\hl{c}{c}) for $p$ does not. Condition~(\hl{b}{b}) then holds for $p{-}1$ by the definition of higher Du Bois singularity, hence condition~(\hl{c}{c}) also holds for $p{-}1$ by inductive hypothesis. We then get the distinguished triangle in $D^b_{\rm coh}(\OO_Y)$
\htt{2.3.1}{}
$$\Gr^F_p\DR_Y^r(\om_Y(*X)/\om_Y)\to\bl(\sigma_{\ges d_Y-p}(\Om_Y^{\ssb}/f\Om_Y^{\ssb},\ddd f\wedge)\br)[d_Y]\to\F\buildrel{+1}\over\to,
\leqno(2.3.1)$$
where $\F$ is a non-zero coherent $\OO_Y$-module supported in $\Sing\,X$. Applying the contravariant functor $\DD$ as in (\hl{1.4.6}{1.4.6}) and using Proposition~\hl{P2.1}{2.1}, we get by an argument similar to 2.1 the distinguished triangle in $D^b_{\rm coh}(\OO_Y)$
\htt{2.3.2}{}
$$(\DD\F)[p{-}d_X]\to\Om_X^p\to\Omu_X^p\buildrel{+1}\over\to.
\leqno(2.3.2)$$
Note that (\hl{2.2.5}{2.2.5}) holds with $p$ replaced by $p{-}2$, hence (\hl{2.2.8}{2.2.8}) holds for $q\eq p$ (although we cannot show the torsion-freeness of $\Om_X^p$). By condition~(\hl{b}{b}) for $p$, there is an isomorphism
\htt{2.3.3}{}
$$\Omu_X^p\cong\Om_X^p.
\leqno(2.3.3)$$
\sk
Set $d_Z:=\dim Z\,$ with $\,Z:={\rm Supp}\,\F$. In order to induce a contradiction, we may assume (by shrinking $Y$ if necessary) that $Z$ is smooth and irreducible, $\F$ is an extension of locally free $\OO_Z$-modules (viewed as $\OO_Y$-modules) so that
\htt{2.3.4}{}
$$\Hc^j\DD\F\eq0\q(j\,{\ne}\,{-}d_Z),\q\h{that is,}\q\DD\F\eq\F\1'[d_Z],
\leqno(2.3.4)$$
where $\F\1'$ is a coherent $\OO_Y$-module supported on $Z$. These imply a distinguished triangle
\htt{2.3.5}{}
$$\F\1'[d_Z{+}p{-}d_X]\to\Om_X^p\to\Om_X^p\buildrel{+1}\over\to.
\leqno(2.3.5)$$
\sk
Set $d_S:=\dim S$ with $S:=\Sing\,X$. Since condition~(\hl{a}{a}) holds for $p{-}1$ (that is, $\alt_X\ges p$) by \cite[Cor.\,1]{hi}, it follows from Proposition~\hl{P1}{1} that
\htt{2.3.6}{}
$$d_Y\mi d_S\ges 2p,\q\h{that is,}\q d_S\pl p\mi d_X\les 1{-}p.
\leqno(2.3.6)$$
By the long exact sequence associated to the distinguished triangle (\hl{2.3.5}{2.3.5}), we get that
\htt{2.3.7}{}
$$d_Z\pl p\mi d_X\eq{-}1\,\,\,\h{or}\,\,\,0,
\leqno(2.3.7)$$
since $\F\1',\Om_X^p$ are $\OO_Y$-modules and $\F\1'\ne 0$. (In these two cases, (\hl{2.3.5}{2.3.5}) induces a short exact sequence.)
\sk
Assume first $d_Z\pl p\mi d_X\eq 0$. Using the long exact sequence of local cohomology sheaves, we get a short exact sequence of coherent $\OO_Y$-modules supported on $Z$\,:
\htt{2.3.8}{}
$$0\to\F\1'\to\Hc_Z^0\Om_X^p\to\Hc_Z^0\Om_X^p\to 0,
\leqno(2.3.8)$$
since $\Hc_Z^0\F\1'\eq\F\1'$, $\Hc_Z^1\F\1'\eq0$ (using for instance \cite[Thm.\,A1.3]{Ei}).
\sk
Let $z$ be the {\it generic point\1} of $Z$ in the scheme-theoretic sense, see \cite{EGA}, \cite{AG}. For a coherent $\OO_Y$-module $\G$ supported on $Z$ in general, we define $\ell_z(\G)$ to be the length of the $\OO_{Y,z}$-module $\G_z$. (This is the same as the sum of the lengths of $\Gr^G_i\G|_U$ with $G$ a finite filtration of $\G$ such that the $\Gr^G_i\G$ are coherent $\OO_Z$-modules (shrinking $Y$ if necessary) and $U\subset Z$ a sufficiently small non-empty open subset such that the $\Gr^G_i\G|_U$ are locally free $\OO_U$-modules.)
\sk
From the short exact sequence (\hl{2.3.8}{2.3.8}) we get that
\htt{2.3.9}{}
$$\ell_z(\F\1')+\ell_z(\Hc_Z^0\Om_X^p)=\ell_z(\Hc_Z^0\Om_X^p),
\leqno(2.3.9)$$
hence $\ell_z(\F\1')=0$, but this is a contradiction.
\sk
Assume now $d_Z\pl p\mi d_X\eq{-}1$. Since $d_Z\les d_S$, we get $p\les 2$ by (\hl{2.3.6}{2.3.6}). Moreover we have
\htt{2.3.10}{}
$$d_S\mi d_Z=0\,\,\,\h{or}\,\,\,1,
\leqno(2.3.10)$$
since (\hl{2.3.6}{2.3.6}) implies that $p\les 0$ if $d_S\mi d_Z\ges 2$.
\sk
In the case $d_S\mi d_Z\eq 1$, we get $p\eq1$ similarly. Here the inequalities in (\hl{2.3.6}{2.3.6}) become equalities so that $d_Y{-}d_S\eq 2$. Then $X$ must be a divisor with normal crossings around a sufficiently general point of $S$, since we have $\alt_X\ges p\eq 1$ by condition~(\hl{a}{a}) for $p{-}1$. In this case, however, $\F$ must be supported on $S$ by Remark~\hl{R1.6}{1.6} (using GAGA). So we get a contradiction.
\sk
In the case $d_S\mi d_Z\eq 0$, we may assume $S\eq Z$. Applying the contravariant functor $\DD$ to (\hl{2.3.5}{2.3.5}), we get the distinguished triangle
\htt{2.3.11}{}
$$\DD\Om_X^p\to\DD\Om_X^p\to\F[d_X{-}p]\buildrel{+1}\over\to.
\leqno(2.3.11)$$
Here $\G^j:=\Hc^j(\DD\Om_X^p)$ ($j\ne-d_X$) is supported in $S$ by Grothendieck duality (applied to the closed immersion $X\,{\setminus}\,S\into Y\,{\setminus}\,S$), see for instance \cite{RD}. Since $p\ges 1$, this implies an exact sequence of coherent $\OO_Y$-modules supported in $Z$\,:
\htt{2.3.12}{}
$$0\to\G^{\1p-d_X}\to\G^{\1p-d_X}\to\F\to\G^{\1p-d_X+1}\to\G^{\1p-d_X+1}\to0,
\leqno(2.3.12)$$
This induces a contradiction using $\ell_z$ as in the case $d_Z\pl p\mi d_X\eq 0$. We can thus proceed by induction on $p$. This finishes the proof of the implication (\hl{b}{b}) $\Longrightarrow$ (\hl{c}{c}).
\par\htt{R2.3}{}\msn
{\bf Remark~2.3.} The definition of higher $p$-log canonicity is independent of the embedding of $X$ into a smooth variety as a hypersurface. This can be verified using the associated analytic space and analytic Hodge ideals \cite{JKSY} together with GAGA. (This may be useful for instance in the toroidal singularity case.) Indeed, if there are two embeddings $X\into Y$, $X\into Y'$, we have locally an isomorphism between $Y^{\rm an}$ and $Y'{}^{\rm an}$ inducing the identity on $X^{\rm an}$. Here we can construct a commutative diagram
$$\begin{array}{ccc}\OO_{Y^{\rm an},x}&\cong&\OO_{Y'{}^{\rm an},x}\\ \downarrow\q\,\,\,\,&&\downarrow\q\,\,\,\,\\ \OO_{X^{\rm an},x}&=&\OO_{X^{\rm an},x}\end{array}$$
by lifting minimal generators of $\mm_{X^{\rm an},x}/\mm_{X^{\rm an},x}^2\,(=\mm_{Y^{\rm an},x}/\mm_{Y^{\rm an},x}^2)$ to $\mm_{Y^{\rm an},x}$ and also to $\mm_{Y'{}^{\rm an},x}$ in such a way that their images in $\mm_{X^{\rm an},x}$ coincide (by lifting them first to $\mm_{X^{\rm an},x}$). Note that morphisms from $\C\{y_1,\dots,y_n\}$ to $\OO_{X^{\rm an},x}$ are determined by the images of the $y_i$.
\bs\bs
\vbox{\centerline{\bf 3. Proofs of Theorems~\hl{T2}{2}--\hl{T3}{3} and Corollary~\hl{C3}{3}}
\bsn
In this section we prove Theorems~\hl{T2}{2}--\hl{T3}{3} and Corollary~\hl{C3}{3} using the nearby and vanishing cycle Hodge modules.}
\msn
{\bf 3.1.~Calculation of a fundamental canonical morphism.} In the notation of 2.1 there is a canonical morphism of filtered right $\D_Y$-modules
\htt{3.1.1}{}
$$({}^r\!M_X,F)\to\bl(\om_Y(*X)/\om_Y,F[-d_X]\br),
\leqno(3.1.1)$$
which underlies the composition of morphisms of mixed Hodge modules
\htt{3.1.2}{}
$$\Q_{h,X}[d_X]\buildrel{\rho}\over\onto{\rm IC}_{h,X}\Q\buildrel{\iota}\over\into\DD(\Q_{h,X}(d_X)[d_X]).
\leqno(3.1.2)$$
Here ${\rm IC}_{h,X}\Q$ denotes the pure Hodge module of weight $d_X$ whose underlying $\Q$-complex is the intersection complex ${\rm IC}_X\Q$. Note that $\iota$ is the dual of $\rho$.
\sk
We have moreover the short exact sequences of mixed Hodge modules
\htt{3.1.3}{}
$$\aligned0\to\Q_{h,X}[d_X]\buildrel{\iota'}\over\longrightarrow\psi_{f,1}\Q_{h,Y}[d_X]&\buildrel{\!\!\rm can}\over\longrightarrow\varphi_{f,1}\Q_{h,Y}[d_X]\to0,\\0\to\varphi_{f,1}\Q_{h,Y}(1)[d_X]\buildrel{\!\!\rm Var}\over\longrightarrow\psi_{f,1}\Q_{h,Y}[d_X]&\buildrel{\rho'}\over\longrightarrow\DD(\Q_{h,X}(d_X)[d_X])\to0,\endaligned
\leqno(3.1.3)$$
which are dual to each other, see \cite[Prop.\,3.8]{mhp} for the surjectivity of ${\rm can}$ and the injectivity of ${\rm Var}$. (These morphisms are defined by $\dd_t$ and $t$ up to a sign at the level of $\D$-modules.)
\par\htt{P3.1}{}\msn
{\bf Proposition~3.1.} {\it We have the equality}
\htt{3.1.4}{}
$$\rho'\ssc\iota'=\iota\ssc\rho:\Q_{h,X}[d_X]\to\DD\bl(\Q_{h,X}(d_X)[d_X]\br).
\leqno(3.1.4)$$
\msn
{\it Proof.} There is a canonical isomorphism (see for instance \cite[(4.5.9)]{mhm})\,:
\htt{3.1.5}{}
$$\Gr^W_{d_X}(\Q_{h,X}[d_X])={\rm IC}_{h,X}\Q.
\leqno(3.1.5)$$
This implies that
\htt{3.1.6}{}
$$\aligned W_{d_X-1}(\Q_{h,X}[d_X])&={\rm Ker}\,\rho,\\ W_{d_X}\DD(\Q_{h,X}(d_X)[d_X])&={\rm Im}\,\iota.\endaligned
\leqno(3.1.6)$$
Recall that $\Q_{h,X}[d_X]$, $\DD\bl(\Q_{h,X}(d_X)[d_X]\br)$ has weights $\ges d_X$ and $\les d_X$ respectively, see \cite[(4.5.2)]{mhm} (which is applied to $a_X:X\to pt$). We then get that
\htt{3.1.7}{}
$$\aligned{\rm Ker}\,\iota'\ssc\rho'&\supset W_{d_X-1}(\Q_{h,X}[d_X]),\\{\rm Im}\,\iota'\ssc\rho'&\subset W_{d_X}\DD(\Q_{h,X}(d_X)[d_X]).\endaligned
\leqno(3.1.7)$$
So the assertion follows, since (\hl{3.1.4}{3.1.4}) holds on the smooth locus of $X$.
\sk
Proposition~\hl{P3.1}{3.1} implies the corresponding assertion in the abelian category constructed in \cite[2.1]{BBD}, which is denoted by $D^b_c(X,\Q)^{[0]}$ in \cite{FPS}.
\sk
We have the following.
\par\htt{T3.1}{}\msn
{\bf Theorem~3.1.} {\it A hypersurface $X\subset Y$ is a $\Q$-homology manifold if and only if the unipotent monodromy part of the vanishing cycle complex $\varphi_{f,1}\Q_Y[d_X]$ vanishes, where $f$ is a function defining locally $X$ in $Y$.}
\msn
{\it Proof.} A complex variety $X$ is called a $\Q$-{\it homology manifold\1} if $H^j_{\{x\}}\Q_X\eq\Q$ for $j\eq 2d_X$, and 0 otherwise for any $x\in X$. It is well-known that this condition is equivalent to that the composition of (\hl{3.1.2}{3.1.2}) is an isomorphism. (This is also equivalent to that the morphism between the underlying $\Q$-complexes is an isomorphism.)
This equivalence can be verified by induction on strata of a Whitney stratification, see also the proof of \cite[Cor.\,1.8]{RSW}.
\sk
The assertion then follows from Proposition~\hl{P3.1}{3.1}, since (\hl{3.1.3}{3.1.3}) implies that
\htt{3.1.8}{}
$$\aligned\Q_X[d_X]&={\rm Ker}\bl(\psi_{f,1}\Q_Y[d_X]\buildrel{\rm can}\over\onto\varphi_{f,1}\Q_Y[d_X]\br),\\\DD(\Q_X(d_X)[d_X])&={\rm Coker}\bl(\varphi_{f,1}\Q_Y(1)[d_X]\buildrel{\rm Var}\over\into\psi_{f,1}\Q_Y[d_X]\br),\endaligned
\leqno(3.1.8)$$
and the morphisms in (\hl{3.1.4}{3.1.4}) is induced by the identity on $\psi_{f,1}\Q_Y[d_X]$. (Indeed, if $\rho\ssc\iota$ is an isomorphism, then the difference between $\bl[\psi_{f,1}\Q_Y[d_X]\br]$  and $\bl[\psi_{f,1}\Q_Y(1)[d_X]\br]$ in the {\it Grothendieck group\1} of mixed Hodge modules vanishes. We then get that $\psi_{f,1}\Q_Y[d_X]\eq 0$ looking at the highest or lowest weight part.) This finishes the proof of Theorem~\hl{T3.1}{3.1}.
\msn
{\bf 3.2.~Condition for bijectivity of graded pieces of (\hl{3.1.1}{3.1.1}).}
Let $({}^r\!M_{\varphi_{f,1}},F)$ be the underlying filtered right $\D_X$-module of $\varphi_{f,1}\Q_{h,Y}[d_X]$. From \cite[(1.3.8)]{hi}, \cite[(2.1.4), (2.2.1)]{mic}, we can deduce the isomorphisms for $p\in\Z$\,:
\htt{3.2.1}{}
$$\Gr^F_{p-d_Y}{}^r\!M_{\varphi_{f,1}}=\Gr^F_{p-d_Y}\Gr_V^0({}^r\B_f,F)=\Gr_{\Vt}^p\OO_Y,
\leqno(3.2.1)$$
where $(^r\B_f,F):=(i_f)_*^{\D}(\Om_Y^{d_Y},F)$ with $i_f:Y\into Y{\times}\C$ is the graph embedding, see \cite{mhp} (and also the proof of \cite[Thm.\.1.5]{MOPW}).
\sk
Replacing $p$ with $p{+}1$ in (\hl{3.2.1}{3.2.1}), we then get the following (since $d_Y\eq d_X{+}1$).
\par\htt{P3.2}{}\msn
{\bf Proposition~3.2.} {\it We have the vanishing}
\htt{3.2.2}{}
$$\Gr^F_{p-d_X}{}^r\!M_{\varphi_{f,1}}=0\q\h{if}\q\alt_f^{\rm min.int}\,{>}\,p{+}1.
\leqno(3.2.2)$$
\sk
Combined with Proposition~\hl{P3.1}{3.1} (and (\hl{3.1.8}{3.1.8})), this gives the following (here the Tate twist~$(1)$ in the source of Var does not cause a problem, since $(F[1])_p\eq F_{p-1}$ by definition).
\par\htt{T3.2}{}\msn
{\bf Theorem~3.2.} {\it If $\alt_f^{\rm min.int}\,{>}\,p{+}1$, there are isomorphism in} $D^b_{\rm coh}(\OO_Y)$\,:
\htt{3.2.3}{}
$$\Gr_F^{d_X-p}\DR_Y^r({}^r\!M_X,F)\simto\Gr^F_p\DR_Y^r(\om_Y(*X)/\om_Y,F).
\leqno(3.2.3)$$
\sk
Note that the left-hand side of (\hl{3.2.3}{3.2.3}) is isomorphic to $\Omu_X^{d_X-p}[p]$ by Theorem~\hl{T1.5}{1.5}.
\msn
{\bf 3.3.~Proofs of Theorems~\hl{T2}{2}--\hl{T3}{3}.} Under the hypothesis of Theorem~\hl{T2}{2}, the right-hand side of (\hl{3.2.3}{3.2.3}) is isomorphic to
$$\bl(\sigma_{\ges d_Y-p}(\Om_Y^{\ssb}|_X,\ddd f\wedge)\br)[d_Y],$$
by \cite[Cor.\,1]{hi} (that is, (\hl{a}{a})\,$\Leftrightarrow$\,(\hl{c}{c}) in Theorem~\hl{T1}{1}), see also (\hl{2.2.1}{2.2.1}) and (\hl{11}{11}). Theorem~\hl{T2}{2} thus follows.
\sk
For Theorem~\hl{T3}{3}, we have the following isomorphism given by Theorem~\hl{T3.2}{3.2}\,:
\htt{3.3.1}{}
$$\Hc^p(\Omu_X^{d_X-p})\cong\Hc^0\Gr^F_p\DR_Y^r(\om_Y(*X)/\om_Y,F).
\leqno(3.3.1)$$
The image of $I_p(X)$ in $\OO_X/(f)$ is isomorphic to a quotient of
$$\Gr^F_p(\om_Y(*X)/\om_Y)\cong I_p(X)/f\1I_{p-1}(X),$$
where $I_{-1}:=\OO_X$. Looking at the differential of the graded de Rham complex, we verify that the image of $I_p(X)$ in $\OO_X/((\dd f),f)$ is isomorphic to a quotient of the right-hand side of (\hl{3.3.1}{3.3.1}). So Theorem~\hl{T3}{3} follows from (\hl{11}{11}) in the introduction.
\par\htt{R3.3}{}\msn
{\bf Remark~3.3.} Assume a hypersurface $X\subset Y$ is a $\Q$-homology manifold, or equivalently, the unipotent monodromy part of the vanishing cycle complex $\varphi_{f,1}\Q_Y[d_X]$ vanishes, see Theorem~\hl{T3.1}{3.1}. Then the morphisms $\iota,\rho$ in (\hl{3.1.2}{3.1.2}) are isomorphisms, hence (\hl{3.2.3}{3.2.3}) is an isomorphism for any $p$. So the assumption $\alt_X\,{>}\,p{+}1$ in Theorem~\hl{T2}{2} can be replaced by $\alt_X\gess p{+}1$, where the latter is needed for the identification of the right-hand side of (\hl{3.2.3}{3.2.3}), see also (\hl{2.2.1}{2.2.1}). (This improvement is, however, rather trivial, since $\alt_X$ cannot be an integer if $\varphi_{f,1}\Q_{h,Y}[d_X]\eq0$.)
\msn
{\bf 3.4.~Spectrum and Tjurina subspectrum.} Assume a hypersurface $X\subset Y$ has only one singular point $x$. The {\it spectrum\1} $\Sp_f(t)=\msum_{i=1}^{\mu_f}\,t^{\1\al_{f,i}}$ can be defined as the Hilbert-Poincar\'e series of the finite-dimensional filtered vector space $\bl(\OO_{Y,x}/(\dd f),V\br)$ so that we have for $\al\in\Q$
\htt{3.4.1}{}
$$\#\{i\in[1,\mu_f]\mid\al_{f,i}\eq\al\}=\dim_{\C}\Gr_V^{\al}\bl(\OO_{Y,x}/(\dd f)\br).
\leqno(3.4.1)$$
Here the filtration $V$ on $\OO_{Y,x}/(\dd f)$ is the quotient filtration of the $V$-filtration on the Brieskorn lattice $H''_f$ (see \cite{Bri}, \cite{Ka2}, \cite{Mal2}, \cite{bl}, etc.), and coincides with the quotient filtration of the microlocal $V$-filtration $\Vt$ on $\OO_{Y,x}$, see \cite[Prop.\,1.4]{JKSY}.
The spectrum $\Sp_f(t)$ depends only on $X$, since it is invariant by a $\mu$-constant deformation as is well-known, see for instance \cite{Va1}.
\sk
The {\it Tjurina subspectrum\1} $\Sp^{\rm Tj}_f(t)\eq\msum_{i\in J_f}\,t^{\al_{f,i}}$ can be defined as the Hilbert-Poincar\'e series of $\bl(\OO_{Y,x}/((\dd f),f),V\br)$ so that
\htt{3.4.2}{}
$$\#\{i\in J_f\mid\al^{\rm Tj}_{f,i}\eq\al\}=\dim_{\C}\Gr_V^{\al}\bl(\OO_{Y,x}/((\dd f),f)\br),
\leqno(3.4.2)$$
where $J_f$ is explained after (\hl{13}{13}) in the introduction.
This subspectrum also depends only on $X$ using an argument similar to \cite[Rem.\,4.2\,(i)]{wh}. (Note that the ideal $((\dd f),f)\subset\OO_Y$ is independent of the choice of $f$.)
\par\htt{R3.4a}{}\msn
{\bf Remark~3.4a.} We have the {\it symmetry} of spectral numbers in the isolated hypersurface singularity case \cite{St1}
\htt{3.4.3}{}
$$\al_{f,i}\pl\al_{f,j}=d_Y\q\h{if}\,\,\,i\pl j\eq\mu_f\pl1,
\leqno(3.4.3)$$
assuming the $\al_{f,i}$ are weakly increasing. It is also well-known that
\htt{3.4.4}{}
$$\alt_X=\min\bl\{\al_{f,i}\mid i\in[1,\mu_f]\br\}.
\leqno(3.4.4)$$
This can be proved by comparing \cite{Mal} and \cite{SS} for instance. (It follows also from (\hl{12}{12}) and \cite[(2.1.4), (2.2.1)]{mic}.)
\par\htt{R3.4b}{}\msn
{\bf Remark~3.4b.} The action of $f$ on $\OO_{Y,x}/(\dd f)$ {\it preserves the filtration $V$ up to the shift by\1} 1, that is,
\htt{3.4.5}{}
$$f\1V^{\al}\bl(\OO_{Y,x}/(\dd f)\br)\subset V^{\al+1}\bl(\OO_{Y,x}/(\dd f)\br)\q(\al\in\Q),
\leqno(3.4.5)$$
since the action of $t$ on the Brieskorn lattice $H''_f$ is defined by that of $f$, see \cite{Bri}, \cite{SS}, \cite[(1.7.2), (2.5.2)]{bl}, \cite[(1.2.6)]{JKSY}, etc.
\sk
From (\hl{3.4.5}{3.4.5}) together with the symmetry (\hl{3.4.3}{3.4.3}) we can deduce the inclusion
\htt{3.4.6}{}
$$f^k\in(\dd f)\q\h{if}\q k\,{>}\,d_Y{-}2\1\alt_X.
\leqno(3.4.6)$$
Indeed, the condition is equivalent to that $\,\alt_X\pl k\,{>}\,d_Y\mi\alt_X$ where the right-hand side is equal to $\max\{\al_{f,i}\}$.
This improves, in the isolated singularity case, the condition $k\gess d_Y$ shown in \cite{BrSk}, see \cite{JKSY2} for details.
\par\htt{R3.4c}{}\msn
{\bf Remark~3.4c.} The {\it semicontinuity\1} of spectral numbers \cite{St2} says that the number of spectral numbers contained in $(a,a{+}1]$ (counted with multiplicity) is {\it upper semicontinuous\1} with respect to a deformation of holomorphic functions with isolated singularities for any $a\in\R$. (Upper semicontinuous means that the number may increase at special points.)
\par\htt{R3.4d}{}\msn
{\bf Remark~3.4d.} Set $\Ht''_f:=\msum_{k\ges 0}\,(\dd_tt)^kH''_f$ (which is called the {\it saturation\1} of $H''_f$). Let $\beta_{f,i}$ ($i\in[1,\mu_f]$) be the eigenvalues of the action of $-\dd_tt$ on $\Ht''_f/t\Ht''_f$. Then $\{-\beta_{f,i}\}$ coincides with the set of roots of $b_f(s)/(s{+}1)$ forgetting the multiplicity, see \cite{Mal}. Moreover, renumbering the $\beta_{f,i}$, we have the well-known {\it integral shifts}\,:
\htt{3.4.7}{}
$$\al_{f,i}-\beta_{f,i}\in\Z_{\ges0}\q(i\in[1,\mu_f]).
\leqno(3.4.7)$$
This can be shown by introducing the filtration $\Fti$ on the Milnor cohomology which is defined by replacing the Brieskorn lattice $H''_f$ with its saturation $\Ht''_f$ in the well-known formula for the Hodge filtration as in \cite{SS}, etc.
\sk
For instance, let $f\eq x^6\pl y^5\pl x^3y^3$. By an argument similar to \cite[(4.1.4)]{bl2}, we see that
\htt{3.4.8}{}
$$\aligned&\{\al_{f,i}\}=\bl\{\tfrac{j}{6}\pl\tfrac{k}{5}\,\big|\,(j,k)\in[1,5]{\times}[1,4]\br\},\\ &\{\beta_{f,i}\}=\{\al_{f,i}\}\cup\bl\{\tfrac{5}{6}\pl\tfrac{4}{5}\mi1,\tfrac{4}{6}\pl\tfrac{4}{5}\mi1\br\}\setminus\bl\{\tfrac{5}{6}\pl\tfrac{4}{5},\tfrac{4}{6}\pl\tfrac{4}{5}\br\}.\endaligned
\leqno(3.4.8)$$
Applying the Thom-Sebastiani type theorem for spectrum and Bernstein-Sato polynomials to $h\eq f\pl g$ with $g\eq z^5\pl w^3$ (see for instance \cite{SS}, \cite[Thm\,0.8]{mic}), we then get that
\htt{3.4.9}{}
$$\alt_h^{\,\rm min.int}=2\,>\,\min\bl(\{\beta_{h,i}\}\cap\Z\br)=1.
\leqno(3.4.9)$$
Indeed, the theorems claim that $\{\al_{h,i}\}=\{\al_{f,j}\pl\al_{g,k}\}$ (similarly for $\{\beta_{h,i}\}$). These imply that $\alt_h^{\,\rm min.int}>\min\bl(\{\beta_{h,i}\}\cap\Z\br)=1$, since
\htt{3.4.10}{}
$$\aligned&\{\al_{f,i}\}\cap\bl[0,\tfrac{7}{15}\br]=\bl\{\tfrac{11}{30}\br\},\q\{\beta_{f,i}\}\cap\bl[0,\tfrac{7}{15}\br]=\bl\{\tfrac{11}{30},\tfrac{14}{30}\br\},\\ &\{\al_{g,i}\}=\{\beta_{g,i}\}=\bl\{\tfrac{j}{5}\pl\tfrac{k}{3}\,\big|\,(j,k)\in[1,4]{\times}[1,2]\br\},\endaligned
\leqno(3.4.10)$$
where $\tfrac{7}{15}=1-\min\{\al_{g,i}\}$. It is easy to see that $\alt_h^{\,\rm min.int}\les 2$.
\msn
{\bf 3.5.~Proof of Corollary~\hl{C3}{3}.} It is enough to show the estimate $\alt_X\less 3/2$ assuming
$$\codim_Y\Sing\,X\eq 4\,\,\,\h{or}\,\,\,5.$$
(Here the codimension is considered in $Y$ instead of $X$.)
Indeed, the upper bound of $\alt_X$ follows from Proposition~\hl{P1}{1} if $\codim_Y\Sing\,X\eq 3$. We can get the lower bound as is remarked after Corollary~\hl{C3}{3}. We may then restrict to an open neighborhood of a sufficiently general point of $\Sing\,X$. We reduce the assertion to the isolated hypersurface singularity case so that Corollary~\hl{C2}{2} can be applied.
\sk
Let $Z$ be a smooth affine variety with $G\subset{\rm Aut}(Z)$ a finite subgroup. Set $X:=Z/G$. Let $G_z\subset G$ be the stabilizer of $z\in Z$. There are Whitney stratification $\Sc,\Sc'$ of $Z,X$ such that $G_z$ is constant on each stratum of $\Sc$, the image of each stratum of $\Sc$ in $X$ is a stratum of $\Sc'$, and $\Sing\,X$ is a union of strata of $\Sc'$. Let $S$ be a maximal-dimensional stratum of $\Sc$ such that its image $S'$ in $X$ is contained in $\Sing\,X$
(Note for instance that the quotient of $\C^n$ by the natural action of the symmetric group ${\mathfrak S}_n$ is smooth.)
Let $z\in S$. Since we are interested in the singularities of $X$ at sufficiently general singular points, we may assume $G\eq G_z$ replacing $G$ if necessary.
\sk
If $\dim S\,{\ne}\,0$, let $h$ be a function such that $h^{-1}(0)$ passes through $z$ transversally to $S$, that is $\ddd h|_S$ does not vanish at $z$. We may assume that $h$ is invariant by $G$ (replacing $h$ with the sum of $g^*h$ for $g\in G$), and moreover $h^{-1}(0)$ is smooth (shrinking $Z$) since $\ddd h|_S\ne 0$ at $z$. Repeating this, we get a transversal space $Z^t$ to $S$ which is invariant by the action of $G$. Divided by $G$, this gives a transversal space $X^t$ to $S'$. Thus the assertion is reduced to the isolated hypersurface singularity case. Hence we may assume
$$\dim\Sing\,X\eq 0,$$
so that Corollary~\hl{C2}{2} can be applied. Note that its first hypothesis $\alt_f^{\,\rm min.int}\,{>}\,p{+}1$ is satisfied in the quotient singularity case, since a $V$-manifold is a $\Q$-homology manifold, that is, 1 is not an eigenvalue of the Milnor monodromy.
\sk
Assume $d_Y\eq 4$ or 5, and
\htt{3.5.1}{}
$$\alt_X>3/2.
\leqno(3.5.1)$$
By Corollary~\hl{C2}{2} for $p\eq 1$, we get
\htt{3.5.2}{}
$$\al_f^{\rm max.Tj}<2,
\leqno(3.5.2)$$
since $\Omu_X^1$ is a sheaf in the quotient singularity case, see \cite[Thm.\,5.3]{DB}.
\sk
From (\hl{3.5.1}{3.5.1}--\hl{3.5.2}{2}) we can deduce that there is no spectral number of $f$ contained in $[2,5/2]$ using Remark~\hl{R3.4b}{3.4b} (applied to $\al\eq\alt_X$). However, this contradicts the {\it symmetry\1} of spectral numbers in Remark~\hl{R3.4a}{3.4a} if $d_Y\eq 4$ (where the center of symmetry is 2).
\sk
If $d_Y\eq 5$, we have no spectral number contained in $[2,3]$ using the symmetry (with center $5/2$). However, this contradicts the {\it semicontinuity\1} of spectral numbers \cite{St2} considering a deformation of $f$ to an $A_1$-singularity $\msum_{i=1}^5\,y_i^2$ (with spectrum $t^{5/2}\1$), see Remark~\hl{R3.4c}{3.4c}. Corollary~\hl{C3}{3} thus follows.
\par\htt{R3.5}{}\msn
{\bf Remark~3.5.} The spectral numbers were once called {\it exponents.} There was, however, a problem of normalization about their range: These are contained either in $(0,d_Y)$ or in $({-}1,d_Y{-}1)$. There is a shift by 1 between the two normalizations. (Note that there is a similar problem about the roots of Bernstein-Sato polynomials.) The former normalization is compatible with jumping coefficients of multiplier ideals, see \cite{Bu}. The latter one comes from the theory of {\it asymptotic\1} mixed Hodge structure \cite{Va1}, and is closely related to the {\it exponents of asymptotic expansions of period integrals.} In order to avoid confusions, we call those contained in $(0,d_Y)$ spectral numbers (although this name is also used in the case of the other normalization by some people).

\end{document}